\documentclass[12pt,leqno,twoside]{amsart}
\usepackage{euscript}
\usepackage{amsmath}
\usepackage{amssymb}

\makeatletter

\renewcommand{\Box}{\diamond}

\theoremstyle{plain}
\newtheorem{Cor[Lem]}{Corollary}

\theoremstyle{definition}

\theoremstyle{remark}

\numberwithin{equation}{section}

\newenvironment{conj}%
{\rm \trivlist \item[\hskip\labelsep{\bf Conjecture}]}%
{\endtrivlist}

\newenvironment{pr}%
{\rm \trivlist \item[\hskip\labelsep{\bf Proof}]}%
{\hspace*{\fill}$\Box$ \endtrivlist}
\newcommand{\prf}{\begin{pr}}
\newcommand{\eprf}{\end{pr}}

\newcommand{\reg}{{\rm{reg}}}

\renewcommand{\int}{{\rm{int}}}

\newcommand{\sing}{{\rm{sing}}}

\newcommand{\cB}{{\EuScript{B}}}

\newcommand{\bR}{{\mathbb R}}
\newcommand{\bC}{{\mathbb C}}

\newcommand{\bZ}{{\mathbb Z}}

\def\t#1{\tilde#1}
\newcommand{\tran}{\bar{\pitchfork}}

\newcommand{\lar}{\longrightarrow}

\newcommand{\ve}{{\varepsilon}}
\newcommand{\vs}{\vskip 3mm}

\begin{document}

\setcounter{section}{0}

\title{Splitting of Singularities}

\author{Guangfeng JIANG$^\dag$  \  {\tiny\rm and } Mihai  TIB\u{A}R}

\thanks{$^\dag$ The research of the first author was  supported by JSPS, NNSFC and BUCT}
\address{G. J.: 
 Department of Mathematics and Information Science,
Faculty of Science,
Beijing University of Chemical Technology,
Bei sanhuan donglu 15,
Beijing 100029,
P. R. China }
\email{jianggf@mailserv.buct.edu.cn}

\address{M.T.: Math\'{e}matiques, UMR 8524 CNRS, Universit\'{e} de Lille 1,
         59655 Villeneuve d'Ascq C\'{e}dex, France }
\email{tibar@agat.univ-lille1.fr}



\begin{abstract} We study one parameter deformations of
a pair consisting of an analytic singular space $X_0$
and a function $f_0$ on it, in case this defines an
isolated singularity. We prove, under general conditions, 
a  bouquet decomposition of the Milnor fibre when the isolated
 singularity splits in the deformation  and the invariance of
 the Milnor fibration if there is no splitting.
\end{abstract}
\maketitle

\section{Introduction}\label{introduction}

Let $f:(X, 0)\lar ({\bC}, 0)$ be an analytic function germ defined on
   an analytic space germ  $(X, 0)$ 
 embedded in $({\bC}^{m+1}, 0)$.  Let
$l$ be a linear function on ${\bC}^{m+1}$, which is  considered  as
the last coordinate function of ${\bC}^{m+1}$.
Let $\Delta$ be  a small open disc in ${\bC}$ with center 0, and
$U$ be a small open  neighborhood of 0 in ${\bC}^{m}$, such that
in $W:=U\times\Delta$,  $(X, 0)$ can be  represented as an analytic set.
For each $t\in \Delta$,
 define $X_t:=W\cap X\cap  l^{-1}(t)$ and
$f_t:=f(-, t)=f|X_t$. Assume that $X_t$ is irreducible  and  $f(0,t)=0$
for any $t\in \Delta$.
 The triple $(X,f,l)$, or briefly, the pair $(X_t, f_t)$, is called
{\em a one-parameter deformation of the (space-function)
pair } $ (X_0, f_0)$.

 Let ${\mathcal S}=\{{\mathcal S}_i\}$
 be a Whitney  stratification of $X$, the representative of $(X,0)$ in $W$.
Denote by $\Sigma _{\mathcal S}(f,l)$ the critical set
 of the mapping $(f,l): X\lar {\bC}^2$
 with respect to the stratification ${\mathcal S}$. We study
deformations of isolated singularities, defined as follows.

\subsection{Definition}\label{definition1}
The triple $(X,f,l)$ (or the pair $(X_t, f_t)$) is called a
{\it  one-parameter deformation of an isolated singularity} 
$(X_0, f_0)$ if the intersection of $\Sigma _{\mathcal S}(f,l)$ with
$l^{-1}(0)$ has   the origin as an isolated point. 

If $(X_t, f_t)$ is a one-parameter deformation of an isolated
singularity   $(X_0, f_0)$, then it follows that  the dimension of
$\Sigma _{\mathcal S}(f,l)$ is at most one and
the intersection of $\Sigma _{\mathcal S}(f,l)$ with $l^{-1}(t)$  is of
dimension 0 (or void) for small $t$.
Since $\Sigma _{\mathcal S}l\subset \Sigma _{\mathcal S}(f,l)$, it
follows that $l^{-1}(t)$ cuts transversally the positive dimensional
strata of ${\mathcal S}$, except at a finite number of points, namely
the points of the set $X_t\cap \Sigma _{\mathcal S}l.$ By the
transversality result of Cheniot \cite{cheniot}, the stratification
${\mathcal S}_t$ of $X_t$ which consists of
 $S_i\cap l^{-1}(t)\setminus \Sigma _{\mathcal S}l$ and the points
$l^{-1}(t)\cap \Sigma _{\mathcal S}l$ is Whitney.
Now  the function
$f_t: X_t \lar {\bC}$ has at most  isolated singularities in $U$
 with respect to the Whitney stratification ${\mathcal S}_t$ of $X_t$.
The critical set of $f_t$ is
$\Sigma _{{\mathcal S}_t}(f_t) = l^{-1}(t)\cap \Sigma _{\mathcal S}(f,l)
$.

\subsection{Definition}\label{splitting}
If $\Sigma _{{\mathcal S}_t}(f_t)$ has only one point in $U$ for small
 enough $U$ and $|t|$,
we say that the singularity $(X_0, f_0)$  {\it does  not split.}

 In this 
case, it follows that the singular locus $\Sigma _{{\mathcal S}}(f,l)$
is non-singular, hence a line up to analytic change of coordinates.

The following questions may
arise in this context.

\begin{conj} {\bf A.}
The Milnor fibre of an isolated singularity $(X_0, f_0)$
is homotopy equivalent to  the bouquet of the Milnor fibres of the isolated
singularities into which it splits.
\end{conj}
\begin{conj} {\bf B.}
If the isolated singularity $(X_0, f_0)$ does not
split, then the Milnor fibration of the isolated
singularity of $(X_t, f_t)$ is homotopically constant,
for $t$ close to 0.
\end{conj}

For  the existence of the Milnor fibrations and the topology of the
Milnor fibre we refer the reader to the papers of L\^{e} \cite{L1,L2}.
There is evidence for these statements as follows. Conjecture A holds when
$X_0$ is a regular space. More generally, it holds when $X_t = X_0$, $\forall
t$, $X_0$ has an isolated
singularity, $\dim X_0\not= 3$ and the singularity of
$f_0$ splits such that outside the origin there are
only Morse singularities, see Siersma's paper \cite{S2}.

We show here that Conjecture A holds in homology (with any
coefficients), for the most general setting. It then
follows, by Whitehead's theorem, that it holds
in homotopy when the singularity splits into only
singular points whose local Milnor fibres are simply
connected.

Conjecture B is an extension of  a well
known result of L\^e and Ramanujam \cite{LR} in case $X=\bC ^m\times
\bC$. The L\^{e}-Ramanujam result has been extended in another direction
by  Vannier \cite{V1, V2} and  Massey \cite{ M1}, \cite{M2}, in case 
$X=\bC ^m\times \bC$  and $f_t$
 with non-isolated singularities on $\bC^{m}$. Let us mention that in
the classical case $X= \bC^m \times \bC$ and $f_t$ with an isolated
singularity, Timourian \cite{Tim} proved furthermore that the
right-equivalence class of $f_t$ is constant. 

With the usual restriction on dimension, we show that Conjectures A and  B
hold  when for each small  $t$ the space $X_t$
has isolated singularity and ``link stability''.

\vs

\noindent{\bf Acknowledgments. } This work started from a visit of the
first author to the
University of Lille 1, and was finished during his stay at Tokyo
Metropolitan University, supported by JSPS. Some discussions with
M. Oka,  B. Teissier and  D. Trotman were helpful for this paper.
He thanks all these institutions and people. The authors thank 
the referee for many valuable remarks that helped to improve the exposition.

\section{Bouquet decompositions}\label{section3}

The purpose of this section is to show that  Conjecture A holds in
homology in general and in homotopy under some conditions. We first extend a 
result of Siersma \cite{S2} on generic splittings to the case of any 
splitting.

\subsection{}\label{2.1}\hskip -2mm
Let $Y$ be an analytic space of pure dimension $n+1$, locally embedded
in ${\bC}^m$ in a neighborhood of $0$. Let ${\mathcal S}$ be a
Whitney stratification
of a representative of $Y$.
Denote by $B$ the open ball in ${\bC}^m$ with radius $\ve $ and   center 0,
by $D$ the  open disk in $\bC$ with  radius $\eta$ and  center 0.

 Let  $f:Y\longrightarrow {\bC}$
be an analytic function  with isolated singularities in  $B\cap Y$
with respect to the stratification ${\mathcal S}$ in the sense of \cite{L2}.
Let  $\Sigma (f)=\{P_0=0, P_1, \ldots , P_{\sigma}\}$ be
 the critical set of $f$ on $B\cap Y\cap f^{-1}(D)$. Denote by
$ b_i=f(P_i)$ ($i=0,1,\ldots, \sigma$).
Assume that $\{P_0=0, P_1, \ldots , P_{\varsigma}\}\subset Y_{\sing} \not= 
\emptyset$,
and $\{ P_{\varsigma+1}, \ldots , P_{\sigma}\}\subset Y_{\reg}$.
Moreover, we assume:

($*$) \quad For all $u\in\bar{D}$, $(Y\cap f^{-1}(u))\bar{\pitchfork}
\partial \bar{B}$ (as stratified sets).

In $B$ (resp. $ D$), take a  small  closed ball $\bar{B}^i$
 $ (\text{resp.  disc }  \bar{D}^i)$ around each $P_i$  (resp. $b_i$)
such that the restriction of $f$ to $\bar{B}^i\cap
 Y\cap f^{-1}(\bar{D}^i\setminus \{b_i\})$ induces the local
Milnor fibration. Let $c_i$ $(1\le i\le \sigma)$ be the path
in  $D\setminus \cup_{j=0}^\sigma\int \bar{D}^j$ connecting $u_0\in \partial\bar{D}^0$ with $u_i\in \partial\bar{D}^i$
 such that each path has no self-intersection and two paths intersect
 only at $u_0$. Without loss of generality, we assume that, if
 $b_i=b_j$ (resp. $b_i=b_0$), then $\bar{D}_i=\bar{D}_j$ and $c_i=c_j$ (resp. $c_i$ is the constant path at $u_0$).
For any $A\subset {\bC}$, denote $Y_A:=Y\cap \bar{B}\cap f^{-1}(A)$. Set
$$
E:=Y _{\bar{D}},\quad \hat{F}:=Y_{u_0}, \quad \quad
   E^i:=\bar{B}^i\cap Y _{\bar{D}^i},\quad
 F^i:=\bar{B}^i\cap Y_{u_i}.
$$

With appropriate deformation retractions and excisions,
 one can prove the following homology
direct sum decomposition formula  which is also true in more
general setting (cf.  \cite{S1, S2, Ti1, J}). 
In this paper we consider  homology with ${\mathbb Z}$-coefficient.

\subsection{Proposition}\label{additivity}
(Additivity of vanishing homology)
{\sl
With the notations and assumptions as above,  we have
$$H_{*}(E, \hat{F})\cong
\bigoplus\limits_{i=0}^{\sigma}
H_{*}(E^i, F^i).
\eqno{\Box}$$
}

\subsection{Decomposition of the fibre in homotopy }\label{2.3}
We make a homotopy model of the wedge of all the local
Milnor fibres $F^i$. Denote
$$
\Gamma=\bigcup\limits_{i=1}^{\sigma}c_i ,\quad
D' =\bigcup\limits_{i=0}^{\sigma} \bar{D}^i ,\quad
E^*=Y_{D'\cup \Gamma} \, \, {\buildrel{h}\over{\simeq}} \, \, Y_{D},
\quad
F^*=Y_{\Gamma}\, \, {\buildrel{h}\over{\simeq}} \, \, Y_{u_0}=\hat{F},
$$
where and in the following, $ \, {\buildrel{h}\over{\simeq}} \,$ means
``is homotopy equivalent to''.

In the fibre $F^*$ one sees the following:

1) $F^0,\ldots, F^{\varsigma}$, the local Milnor fibres of $f$ at $P_0,
\ldots , P_ {\varsigma}$;

2)  the vanishing cycles  from each
$F^{\varsigma +j}\, \, {\buildrel{h}\over{\simeq}} \, \,
S^n\vee \cdots \vee S^n$
 ($\beta _{\varsigma +j}$ copies of $n$-sphere), the
local Milnor fibre of $f$ at each $P_{\varsigma +j}\in Y_{\reg}$,
where $\beta _{\varsigma +j}$ is the local Milnor number of $f$ at
$P_{\varsigma +j}$.

 Let $h_1^{\varsigma +j}\cup \cdots \cup
 h_{\beta _{\varsigma +j}}^{\varsigma +j}$ 
be the $(n+1)$-cells (called the thimbles) to
be attached to $F^{\varsigma +j}$ in order to kill the vanishing cycles.
Let $H$ be the union of all the thimbles over all $j$.

Assume that $\hat{F}$ and $F^i$ ($0\le i\le \varsigma$)
 are connected. Let $x_i\in \partial F^i$ ($0\le i\le \varsigma$),
 and let $x_{\varsigma +j}\in F^{\varsigma +j }$
be the wedge point of the spheres. Take a non self-intersecting path $\gamma _i$ in $F^*$ connecting
$x_0$ and $x_i$ $( 1\le i\le \sigma )$ by lifting $c_i$ (if $b_i \not= b_0$) 
or within $Y_{u_0} = \hat F$ (if $b_i = b_0$), so that
two paths intersect only at $x_0$. We also want that $\gamma _i$ does not intersect $F^j$, for $j\not= i$. In order to satisfy this condition, we may need to modify the path within a tubular neighbourhood of $F^*$, resp. $\hat F$, which is of course possible. We then have the inclusion
$$\iota :F':=F^0\cup(\gamma _1\cup F^{1} )\cup \cdots
 \cup (\gamma _{\sigma}\cup F^{\sigma} )
\hookrightarrow F^*.$$
Note that
$F'\,\,{\buildrel{h}\over{\simeq}}\,\, F^+\vee S$, where
$$F^+:=
F^0\cup(\gamma _1\cup F^{1} )\cup \cdots\cup (\gamma _{\varsigma}\cup
 F^{\varsigma})
 \, \,{\buildrel{h}\over{\simeq}} \, \,
F^0\vee F^1\vee \cdots \vee F^{\varsigma},$$  and
$S:=S^n_1\vee \cdots \vee S^n_{\beta} $ is the wedge of
 $\beta=\sum\limits_{j} \beta_{\varsigma +j}$
copies of the $n$-sphere.

 Let $B'_j$ be the ball with boundary $S^n_j$ in the bouquet of spheres.
 Then we have the  inclusion
$$F'\hookrightarrow  F^+\vee B':= F^+\vee B'_1\vee \cdots \vee B'_{\beta}.$$
Define
$\varphi: F^+\vee S \hookrightarrow F^*$ by the composition of
$F^+\vee S \,\,{\buildrel{h}\over{\simeq}} \,\, F'
{\buildrel{\iota}\over{\hookrightarrow}}F^*$. 
 From the identification of balls with  thimbles
we obtain the following maps
$$
\varphi': F^+\vee B' \lar F^*\cup H, \quad
\varphi'':F^+   \hookrightarrow F^+\vee B' \lar F^*\cup H. \quad $$

\subsection{Theorem}\label{bouquet1} 
{\sl Under the above assumptions,
if $E$ is contractible, then
the map $\varphi$ induces isomorphisms on all
 the homology groups
$$H_*(F^+\vee S){\cong} H_*(F^*).$$
Moreover, if $F^0, \ldots, F^{\varsigma}$ and $\hat{F}$ are simply connected,
then
$$\hat{F}\,\,{\buildrel{h}\over{\simeq}} \,\,F^+\vee S
\,\,{\buildrel{h}\over{\simeq}} \,\,F^0\vee \cdots \vee F^{\varsigma}
 \vee S. $$
}

\prf
 We follow the proof of \cite[Proposition 2.8]{S2}.
The  maps $\varphi $ and $\varphi '$ above give a map  between the space
pairs:
$$\varphi ^{\mathrm rel}: (F^+\vee B', F^+\vee S)\lar
(F^*\cup H, F^*),$$  which induces the map between the homology groups:

\vs
{\setlength{\unitlength}{1pt}

\begin{picture}(380,130)(-120,-50)
\put(-115, 62){$H_{q+1}(F^+\vee B', F^+\vee S)$}
\put(10,65){\vector(1,0){40}}
\put(65,62){$H_q( F^+\vee S )$}
\put(135,65){\vector(1,0){40}}
\put(188,62){$H_q(F^+\vee B' )$ }
\put(-50,50){\vector(0,-1){30}}\put(-40,35){$\varphi ^{\mathrm rel}_*$}
\put(80,50){\vector(0,-1){30}}\put(82,35){$\varphi _*$}
\put(205,50){\vector(0,-1){30}}\put(207,35){$\varphi '_*$}
\put(-95, 2){$H_{q+1}(F^*\cup H, F^*)$}
\put(10,5){\vector(1,0){40}}
\put(68,2){$H_q(F^*)$}
\put(135,5){\vector(1,0){40}}
\put(183,2){$H_q(F^*\cup H  ).  $ }
\put(55,-35){Diagram 1}
\end{picture}
}

By excision, $\varphi ^{\mathrm rel}_*$ is an isomorphism
(cf. \cite[\S3]{La}).
Note that  $F^+ \,\,{\buildrel{h}\over{\simeq}} \,\, F^+\vee B'$.
By mainly excisions, it follows that  the inclusion
$(E^+, F^+)\hookrightarrow (E^*, F^*\cup H)$ induces an isomorphism in
homology, where
$$E^+: =\left(\bar{B}^0\cap Y\cap f^{-1}(\bar D^0)\right)\cup
\left(\bigcup\limits_{i=1}^{\varsigma}\gamma _i\right)\cup
\left(\bigcup\limits_{i=1}^{\varsigma}
 \bar{B}^i \cap Y\cap f^{-1}(\bar{D}^i)\right).$$ Hence
 $\varphi '_*$ is an isomorphism since both $E^+$ and $E^*$
are contractible. These imply that $\varphi _*$ is an isomorphism.
\eprf

\vs

We return to our original settings.
Let $(X,0)$ be an analytic space germ of dimension $n+1>2$,
 locally embedded in $({\bC}^{m+1},0)$.  Let $f:(X,0)\longrightarrow
({\mathbb C},0)$ be a function germ. 
Let $l$ be a linear function,
considered as the last coordinate of ${\bC}^{m+1}$,
and denote $X_t=X\cap l^{-1}(t)$.
 The definition (Definition~\ref{definition1}) 
of  one-parameter deformation  $(X_t, f_t)$ of an isolated
singularity $(X_0, f_0)$ implies the following facts:
\begin{itemize}
\item[(1)] $l^{-1}(0) $ intersects all the strata of $X\setminus \{0\}$
transversally. Note that
 the strata of dimensions less than 2 are contained in 
$\Sigma_{\mathcal S}(f,l)$. For any stratum $S_i\in {\mathcal S}$ of
 dimension at least 2 and any point $z\in S_i\cap l^{-1}(0)$,
if the transversality fails, then $z$ is  a critical point
of $l|_{S_i}$, the restriction of $l$ to $S_i$. Since
the critical locus $\Sigma l|_{S_i}$ of $l|_{S_i}$ 
is contained in $\Sigma (f,l)|_{S_i}\subset \Sigma_{\mathcal S}(f,l)$ and
$l^{-1}(0)\cap \Sigma_{\mathcal S}(f,l)=\{0\}$, we have $z=0$;
\item[(2)]  $l^{-1}(0)$ is transversal to all the strata of
	   $f^{-1}(0)\setminus \{0\}$.
\end{itemize}

By applying Proposition \ref{additivity} to  $(X_t, f_t)$, we see
immediately from the following Lemma \ref{basic lemma} that Conjecture A
is true in homology.

\subsection{Lemma}\label{basic lemma}
{\sl
 Let $(X_t, f_t)$ be a one-parameter deformation of
 the  isolated singularity $(X_0, f_0)$. Then we have
\begin{itemize}
\item[1)]{ For any $\ve >0$ small enough, $f^{-1}_0(0)\cap X_0$ 
	   is transversal to  the boundary $\partial \bar{B}_{\ve}$ 
	   of the closed ball $\bar{B}_{\ve}\subset \bC^m$ with center $0$ and 
	   radius $\ve $. We denote this by 
           $(f^{-1}_0(0)\cap X_0)\tran \partial \bar{B}_{\ve}$;}
\item[2)]{Fix an $\ve _0$ with the property in 1). There exist
	$\eta>0,\tau>0$ such that for any $|u|<\eta$ and $|t|<\tau$, we
	have $(f^{-1}_t(u)\cap  X_t)\tran\partial\bar{B}_0$, where
$\bar{B}_0:=\bar{B}_{\ve_0}$;}

\item[3)]{Let $\ve _0>0$, $\eta _0>0$ and $\tau_0>0$ be as in 2).
         If $\tau_0$ is small enough, then for any $|t|<\tau _0$ and 
	 $u\in \partial \bar{D}_{0}$,
         $\hat{F}:=f^{-1}_t(u)\cap X_t\cap \bar{B}_{0}
          \, \, {\buildrel{h}\over{\simeq }} \, \,
F_0:=f^{-1}_0(u)\cap X_0\cap \bar{B}_{0}$, where $\bar{D}_0$ is the 
	closure of the open disk in ${\mathbb C}$ with center $0$ and
	radius $\eta _0$;}
\item[4)]{$f_t^{-1}(\bar{D}_{0})\cap X_t\cap \bar{B}_{0}
\, \, {\buildrel{h}\over{\simeq }} \, \,
 f_0^{-1}(\bar{D}_{0})\cap X_0\cap \bar{B}_{0}$. In
	particular, if $\eta _0$ is small enough, both spaces are 
         contractible.}
\end{itemize}

}

\begin{pr} Part 1) is the well known lemma of the ``conic structure of the 
 analytic germs'', see, for instance, \cite{Milnor} for the smooth case
 and  \cite{bv} for the stratified case.

The proof of statement 2)  follows from \cite[\S2]{L1}.
Note that the conditions required for $l$ in loc. cit. such that the
proof works are fulfilled  by our $l$,
 since $\dim \Sigma_{\mathcal S}(f,l)\le 1$
(see also the similar remarks in \cite[\S1.1]{Ti1}).

To prove 3), we consider the map (cf. \cite{L1,L2})
$$G=:(f, l):X\cap(\bar{B}_{0}\times \Delta)\lar{\bC}\times \Delta,$$
where $\Delta $ is the open disc in ${\bC}$ with center $0$ and radius
$\tau _0$.
Let
$$Z_i^{(1)}=G^{-1}(\partial \bar{D}_{0}\times \Delta)
\cap {\mathcal S}_i\cap (B_{0}\times \Delta),
\, \,  \text{ and } \, \,
Z_i^{(2)}=G^{-1}(\partial \bar{D}_{0}\times \Delta)\cap {\mathcal S}_i\cap (
\partial\bar{B}_{0}\times \Delta)$$
be the strata of  the Whitney stratification of
$G^{-1}(\partial \bar{D}_{0}\times \Delta)\cap X\cap
 (\bar{B}_{0}\times \Delta)$
induced  from ${\mathcal S}=\{{\mathcal S}_i\}_i$. Obviously,
each $G|Z_i^{(1)}$ is a submersion. By the transversality 2),
each  $G|Z_i^{(2)}$ is  again a  submersion.
 It follows from  Thom-Mather's first isotopy lemma that 3) holds.

Let $Z:=G^{-1}(\bar{D}_{0}\times \Delta)\cap X\cap
(\bar{B}_0\times \Delta)\lar \Delta$
be the map $\pi \circ G$,
where $ \pi: \bar{D}_{0}\times \Delta\lar \Delta $ is the
projection to the second component. Stratify $Z$ by
$Z_i^{(1)}, Z_i^{(2)} $ and
$$
Z^{(3)}_i=  G^{-1}(D_{0}\times \Delta)
\cap {\mathcal S}_i\cap ( B_0\times \Delta ), \quad
Z^{(4)}_i = G^{-1}(D_{0}\times \Delta)\cap {\mathcal S}_i\cap
(\partial\bar{B}_0\times \Delta ).
$$
It is clear that  $\{Z^{(j)}_i\}$ is  a Whitney stratification of
 $Z$ and the restrictions of $\pi$ to each stratum is a submersion.
By Thom-Mather's first isotopy lemma,
$\pi$ is a locally trivial topological  fibration.
This proves 4).
\end{pr}

\vs

In some cases, one can remove  from
Theorem \ref{bouquet1} the requirement that  the Milnor fibres $F^i$
  be  simply connected, and get a bouquet
decomposition in homotopy. For example, if one can prove that
the map $\varphi$ induces isomorphisms on the fundamental groups of the spaces
and on the homologies of the universal coverings of the spaces, then use
Whitehead's theorem \cite{Wh}. This is the approach of Siersma \cite{S2}.
In the remainder of this section we use this idea to prove that
under some assumptions Conjecture A is also true in homotopy.

Define $\rho: {\bC}^{m+1}\lar {\bR}$ by
$\rho (z_1, \ldots, z_m, z_{m+1}):=\sum\limits_{j=1}^{m} z_j\bar{z_j}$.
Still denote by $\rho$ its  restriction to $X$.
 Denote by $\Gamma _{\mathcal S}(\rho, l)$ the germ of the set
$\overline{\Sigma_{\mathcal S}(\rho,l)\setminus \rho ^{-1}(0)}$ at the
origin,  where $\Sigma_{\mathcal S}(\rho,l)$ denotes the critical set of 
the map $(\rho,l): X\longrightarrow {\bR\times \bC}$ relative to the
 stratification ${\mathcal S}$ of $X$.

\subsection{Definition}\label{linkstability} Let $X_t$  be a space
family with $0\in X_t$ for all $t\in {\mathbb C}$. Identify
 ${\rm Cone}\,(X_t\cap \partial B_0)$ with
$(X_t\cap \partial B_0)\times [0,\ve _0]/(x,0)\sim (y,0)$,
 where $ B_0$ is the open ball with center 0
and radius $\ve _0$.
If there exist $\ve_0>0$ and $\tau _0>0$ such that for each $|t|<\tau_0$,
there exists a homeomorphism $\varkappa _t$ from
 ${\rm Cone}\,(X_t\cap \partial \bar{B}_0)$   to $X_t\cap \bar{B}_0$ such that 
$\rho\circ \varkappa _t$ is the
projection onto the interval $[0,\ve _0]$, we say that the family of
 germs $(X_t, 0)$
has {\it link stability}.

\subsection{Lemma}\label{cone}
{\sl
If $\Gamma _{\mathcal S}(\rho ,l)= \emptyset$,
then  the family $(X_t,0)$ has link stability.
}

\prf
Let $W$ be an open neighborhood of the origin of ${\bC}^{m+1} $ such
that inside $W\cap X$,   $\Gamma _{\mathcal S}(\rho ,l)= \emptyset$.
There exist $\ve _0>0, \tau _0>0$
such that $\bar{B}_0 \times \Delta \subset W$,
 where $B_0$ is the open  ball in
${\bC}^m $ with center 0 and radius $\ve _0$, and $\Delta$ is the open disc
in ${\bC}$ with center 0 and radius $\tau _0$. Since
$\Gamma _{\mathcal S}(\rho ,l)= \emptyset$,
for each $t\in \Delta$, the restriction of $\rho$
to each stratum of $X_t$  is a submersion, except at the origin.
 In other words, there are no 0-dimensional strata of $X_t$ except
the origin $(0,t)$ and  each positive dimensional stratum
of $X_t$ intersects $\partial\bar{B_{\ve}}$
transversally, for each $0<\ve \le\ve _0$.



By \cite[II (3.3)]{gibson} or \cite[p. 42]{Goresky-Mac},
 there exists a controlled vector field $v$ on a
 punctured neighborhood  $U\setminus \{0\}$ of $\bar{B}_0\setminus
 \{0\}$,
tangent to the strata of $X_t$,
 such that
$$d_z\rho(v)=-\left(\frac{d}{ds}\right)_{\rho},$$
where $U$ is an open neighborhood of $\bar{B}_0$ and  
$\left(\frac{d}{ds}\right)_{\rho}$ is the unit tangent  vector
to  $\bR$ at $\rho$.

By \cite[II (4.7)]{gibson} or \cite[p. 42]{Goresky-Mac}, this vector 
field $v$ can be integrated, and by choosing the
initial values appropriately, we can get the desired homeomorphism.

More precisely, let $y\in X_t\cap \partial\bar{B}_0$, and
$h_y:(-\delta,\delta)\lar X_t\setminus \{0\}$ be the integral
curve of $v$ with
$h_y(0) =y$. The  following points  are  important  for the
integration of the controlled vector fields.
For each $s\in (-\delta,\delta)$, $h_y(s)$ is the unique point
on the orbit passing through $y$, and $h_y(-\delta,\delta)$ is in the
same stratum which contains $y$.

On each stratum, $h_y$ is smooth, and $\frac{d(\rho\circ h_y)}{ds}(s)=-1$, so
$ \rho\circ h_y(s)=\rho(y) -s=\ve _0-s$. Define
$$\varkappa _t :(0, \ve _0]\times (X_t\cap\partial\bar{B}_0)\lar
(X_t\cap \bar{B}_0)\setminus\{0\}, \quad \mathrm{ by } \quad
(s, y)\longmapsto h_y(\ve _0-s).$$
Then $\rho\circ \varkappa _t (s, y)=s$.
This also shows that,  $\varkappa _t $ is smooth on each stratum, and
$\rho\circ \varkappa _t$ is the projection onto $(0, \ve _0]$.
Hence $\varkappa _t $ can be extended to a
 homeomorphism between ${\rm Cone} (X_t\cap \partial\bar{B}_0)$ and
$ X_t\cap \bar{B}_0$.
\eprf

\subsection{Remark}\label{triviality}
 The assumption
 $\Gamma _{\mathcal S}(\rho ,l)= \emptyset$
in Lemma \ref{cone} is satisfied  automatically in some
cases: $X=X_0 \times {\bC}$, or $X_t$ is a family of weighted
homogeneous complete intersections with isolated singularity.
If $X_t$ is a family of hypersurfaces with isolated singularities and
does not split, then
this assumption implies the topological triviality of the family.
 This is similar to the so-called ($m$)
condition for the pair $(X\setminus (0\times \bC), 0\times \bC )$
used  in the literature
(see, {\it e.g.}, \cite{BK}). However, this condition
does not imply that the pair $(X\setminus (0\times {\bC}), 0\times {\bC})$
satisfies Whitney condition  as shown by the following example.

\subsection{Example}\label{example4-3} 
Let  $X$ be the Brian\c{c}on-Speder family of surfaces 
defined by
$h=z^3+ty^{2\alpha +1}z+xy^{3\alpha +1}+x^{6\alpha +3}=0$ $(\alpha\ge 1)$.
Then $(X\setminus (0\times {\bC}), 0\times {\bC})$ does not
satisfy the Whitney  condition (cf. \cite{BS}).
Computation shows that  $\Gamma _{\mathcal S}(\rho ,l)= \emptyset$
(see also Example \ref{5.3}).

\vskip 3mm

\subsection{Theorem}\label{bouquet22}
{\sl
 Let $(X_t, f_t) $ be a one-parameter  deformation of the  isolated
singularity $(X_0, f_0) $, with $X_t$ irreducible at $0$ and $\dim X_t=n+1\ne 
3$, $\forall t$, $n\ge 1$.
 Suppose that there exist $\ve _0>0$ and $\tau _0>0$ such that 
$X_t\cap B_0 \setminus \{0\}$ is non-singular for all $t\in \Delta$,
and that $l^{-1}(t)\cap \Sigma _{\mathcal S}(f,l)
:=\{P_0(t)=0, P_1(t),\ldots, P_{\sigma}(t)\}$ for $t\in
 \Delta\setminus \{0\}$ and $\lim\limits_{t\rightarrow 0}P_j(t)=0$.
If  $(X_t,0)$ has link stability, then
$$F_0 \, \,{\buildrel{h}\over{\simeq}}\,\, F_t\vee S^n\vee\cdots
\vee S^n,$$
where $F_t$ is the Milnor fibre of $f_t$ at $0$, and
 the total  number of spheres $S^n$ in the bouquet is equal to 
the sum of the local Milnor numbers of $f_t$ at $P_i(t)$.
}

\prf  We  use Theorem~\ref{bouquet1} and
 Lemma \ref{basic
lemma}. Let $\ve _t >0$ and $\eta _t>0$ be the Milnor data for $f_t$,
i.e., the restriction 
$f_t : \bar{B}_t\cap X_t\cap f_t^{-1}(\bar{D}_t^*)\lar
\bar{D}_t^*:=\bar{D}_t \setminus {0}$
of $f_t$ is the Milnor fibration of $f_t$,
where $B_t$ is an  open ball with center $0$ and radius $\ve _t$, and 
$D_t$ is an open disc with center $0$ and radius $\eta _t$. We also
assume that $0<\ve _t<\ve _0$ and $0<\eta _t<\eta _0$ for $t\ne 0$. 

  We briefly recall the constructions in \S\ref{2.1}
and \S\ref{2.3} for $(X_t, f_t)$.
By the assumptions,
there exists  
$\tau _0>0$ such that for any $t \in \Delta $, $X_t$ has an isolated 
singularity in 
${\cB}:=X_t\cap \bar {B}_0$.

Set  $b_0:=f_t(P_0(t))=0$, 
$b_i:=f_t(P_i(t)) \in D_0$.
Note that $\varsigma=0$, since $X_t$ has an isolated singularity at
$P_0(t)=0 $ in  ${\cB}$. 
In $D_0$, take small closed discs $\bar{D}^i$ with center $b_i$ and
radius $\eta '>0$. Let $u_i$ be a point on $\partial\bar{D}^i$. For
$i>0$, let
 $c_i$ be the path connecting $u_0$ with $u_i$, as explained in \S\ref{2.3}.

 Let
$\bar{B}^i$ be the closed ball with center $P_i(t)$  and radius $\ve'>0$.
Take $\eta'>0, \ve '>0$  so small  that 
$\bar{B}^i$ $ (\text{resp. }  \bar{D}^i)$  are disjoint and
contained in $B_0$ (resp. $ D_0$),
the restriction of $f_t$ to 
$X_t\cap\bar{B}^i \cap f_t^{-1}(\bar{D}^i\setminus\{b_i\})$
is a Milnor fibration with fibre $F^i:=X_t\cap \bar{B}^i\cap f_t^{-1}(u_i)$,
 and $ E^i:=\cB\cap\bar{B}^i\cap f^{-1}(\bar{D}^i)$
is contractible.

Similarly, one has $E$, $\hat{F}$, $E^*$,
 $F^* \, \,{\buildrel{h}\over{\simeq}}\,\, \hat{F} $,
 $F^+=F^0$, and the mapping
$\varphi: F^0\vee S\hookrightarrow F^*$. Note that $F^i$ is connected, since 
$X_t$ is irreducible at $0$ and with isolated singularity, $\forall t$.

Since  $X_t\cap B_0$ has an isolated singularity at $0$,
 $F^+=F^0$ is the Milnor fibre $F_t$ of $f_t$. 
If $\dim X_t=2$, by using resolution of singularities, one can prove
that  $F_t$ is a bouquet of one-spheres. Hence, the theorem follows.
If $\dim X_t>3$ and  $F_t$ is simply connected, then the theorem
  also follows from  Theorem \ref{bouquet1} and Lemma \ref{basic lemma}.
  
In the general case, we make use of Whitehead's theorem \cite{Wh} in a similar manner as done by Siersma in \cite{S2}.
Namely, we prove the following two statements in the remainder of this section:
\begin{itemize}
\item[1)]{the map $\varphi$ induces isomorphism on the
fundamental groups
(cf. Lemma \ref{lemma5});}
\item[2)]{the  map $\varphi$ induces isomorphism on the homology groups of
universal coverings of the spaces (cf. Proposition  \ref{covering}).}
\end{itemize}

Then we may apply Whitehead's theorem \cite{Wh} to conclude that $\varphi$
is a homotopy equivalence.
\eprf

\subsection{Lemma}\label{lemma5}
{\sl
 Let $(X_t, f_t) $ be a one-parameter deformation of the  isolated
singularity $(X_0, f_0) $ with  $\dim X_t> 3$.
If  $X_t\cap B_0$ has an isolated
 singularity at $0$ and  $(X_t,0)$ has link stability, then
\begin{itemize}
\item[1)]{$\pi _1(\partial F_t)\cong
\pi _1(\partial F^*)$, where $\partial F^*:=F^*\cap \partial \bar{B}_0$;}
\item[2)]{$\pi _1( F_t)\cong \pi _1(F^*)$, hence
 $\pi _1( F_t\vee S)\cong \pi _1(F^*)$. }
\end{itemize}
}

\prf
 Since $\hat{F}$ is a deformation retract of $F^*$
(cf. \S\ref{2.3}), it is enough to prove the lemma by replacing
 $F^*$ by $\hat{F}$; i.e.
\begin{itemize}
\item[$1'$)]{ $\pi _1(\partial F_t)\cong
\pi _1(\partial \hat{F})$;}
\item[$2'$)]{
$\pi _1( F_t)\cong \pi _1(\hat{F})$, hence
 $\pi _1( F_t\vee S)\cong \pi _1(\hat{F})$.}
\end{itemize}

We follow \cite[\S3]{S2} closely.
In the proof, we use the following notations:
$$\begin{array}{lll}
M=\partial \bar{B}_0\cap X_t, \,\,&
 \hat{F}=f_t^{-1}(u_0)\cap\bar{B}_0\cap X_t,
\,\,& K=f_t^{-1}(0)\cap\partial\bar{B}_0\cap X_t,\cr
M_t=\partial\bar{B}_t \cap X_t, \,\,&
 F_t=f_t^{-1}(u_0)\cap\bar{B}_t\cap X_t ,
\,\,& K_t=f_t^{-1}(0)\cap\partial\bar{B}_t\cap X_t,\cr
\triangledown F=\overline{\hat{F}\setminus F_t} ,\,\,&
\triangledown M=\overline{(B_0\setminus B_t)\cap X_t}. \,\,
\end{array}
$$
Note that we have
$ K \, \,{\buildrel{h}\over{\simeq}}\,\,\partial \hat{F}$ and 
$  K_t \, \,{\buildrel{h}\over{\simeq}}\,\,\partial F_t$.

{\setlength{\unitlength}{1pt}

\begin{picture}(380,170)(-120,-90)

\put(-30, 62){$\pi _1(M_t)$}
\put(15,65){\vector(1,0){40}}\put(30,67){$\psi _1$}
\put(65,62){$\pi _1(\triangledown M)$}
\put(152,65){\vector(-1,0){40}}\put(127,67){$\psi _2$}
\put(163,62){$\pi _1(M)$ }
\put(-17,15){\vector(0,1){40}}\put(-15,30){$\cong$}
\put(80,15){\vector(0,1){40}}\put(82,30){$\phi _1$}
\put(175,15){\vector(0,1){40}}\put(177,30){$\cong$}
\put(-30, 2){$\pi _1(K_t)$}
\put(15,5){\vector(1,0){40}}\put(30,10){$\phi _2$}
\put(65,2){$\pi _1(\triangledown F)$}
\put(152,5){\vector(-1,0){40}}\put(127,7){$\cong$}
\put(163,2){$\pi _1(K)$ }
\put(-17,-5){\vector(0,-1){40}}\put(-15,-30){$\cong$}
\put(80,-5){\vector(0,-1){40}}\put(82,-30){$\cong$}
\put(-30, -58){$\pi _1(F_t)$}
\put(15,-55){\vector(1,0){40}}\put(30,-50){$\phi _3$}
\put(65,-58){$\pi _1(\hat{F})$}

\put(55,-86){Diagram 2}
\end{picture}
}

All the morphisms in  Diagram 2 are induced by the inclusion maps.
The indicated isomorphisms can be proved via Morse theory, by using the results of Hamm \cite[2.9]{Ha}.
By link stability, the inclusions of  $M_t$ and $M$ into  $\triangledown M$
 are homotopy equivalences. We have the isomorphisms $\psi _1$ and $\psi _2$.
It follows that $\phi _1$ and $\phi _2$ are   isomorphisms. Hence
  $\phi _3$ is an isomorphism.
\eprf

\subsection{}\label{4.4}\hskip -2mm
 Let $(X_t, f_t) $ be a one-parameter  deformation of the  isolated
singularity $(X_0, f_0) $ with  $\dim X_t> 3$. Assume
$X_t\cap B_0$ has an isolated singularity at $0$  and $(X_t,0)$ has 
link stability.
We continue to use the notations in \S\ref{bouquet22} and \S\ref{lemma5}.

 By link stability, the  cone $cM$ over
$M$ is homeomorphic to ${\cB}$. Let
 $\t{M}$ be  the universal covering of $M$. $\t{M}$  is smooth,
connected,  and  simply connected.
 Set $\t{{\cB}}:=c\t{M}$, the cone over $\t{M}$,
which is smooth outside the top $*$. There is a map
$\pi : \t{{\cB}}\lar {\cB}$  compatible with the cone structure such that
the restriction of $\pi$ to $\t{{\cB}}\setminus {*}$ is also a covering, which
can be identified with $(0, 1]\times \t{M} \lar (0, 1]\times M$. The
function $f_t$  on $\bar{B}_0\cap X_t$ and its
 restriction to $M$ can be lifted to
functions on $\t{{\cB}}$ and $\t{M}$ respectively; i.e., we have commutative 
Diagram 3.

{\setlength{\unitlength}{1pt}

\begin{picture}(380,120)(-120,-35)
\put(-15, 62){$\t{M}$}
\put(20,65){$\hookrightarrow$}
\put(65,62){$\t{{\cB}}$}
\put(95,60){\vector(1,-1){40}}\put(125,40){$\t{f}_t:=f_t\circ \pi$}

\put(-7,55){\vector(0,-1){40}}\put(-5,30){$\pi$}
\put(70,55){\vector(0,-1){40}}\put(77,30){$\pi$}

\put(-15, 2){$M$}
\put(20,5){$\hookrightarrow$}
\put(65,2){${\cB}$}
\put(90,5){\vector(1,0){40}}\put(110,9){$f_t$}
\put(150,2){$\bC$ }
\put(55,-30){Diagram 3}
\end{picture}
}

\subsection{Lemma}\label{lemma4.4}
{\sl
Under the assumptions above,  $\t{K}=\t{f}_t^{-1}(0)\cap \t{M}$ and 
$ \tilde{\hat{F}}=\t{f}_t^{-1}(u_0)\cap \t{{\cB}}$ are simply
connected. Moreover, the restrictions of $\pi$ give universal coverings
$$\pi : \t{K}\lar K  \quad \text{ and } \quad\pi: \tilde{\hat{F}}\lar 
\hat{F}.$$
}

\prf
  One  uses the following fact from
 topology (see \cite[Lemma 4.1]{S2}).

\noindent{\bf Sublemma } Let $Y$ be connected, locally path connected
and locally simply connected. Let $\pi :\t{Y}\lar Y$ be the  universal
covering of $Y$. If $Z(\subset Y)$ is connected and the inclusion map
$Z\hookrightarrow Y$ induces an isomorphism on the fundamental groups,
then the restriction of $\pi$ to $\t{Z}:=\pi ^{-1}(Z)$ is also a universal
covering.

The space ${\cB}\setminus\{0\}\,\, {\buildrel{h}\over{\simeq}} \,\, M$
 satisfies the requirements for $Y$ in the sublemma. By Lemma \ref{lemma5},
 both $K$ and $ \hat{F}$
satisfy the requirements for $Z$ in the sublemma, and  the lemma follows.
\eprf

\subsection{}\label{4.5}\hskip -2mm
  We  repeat the constructions in \S\ref{2.1}
 and \S\ref{2.3}
for the spaces $\t{{\cB}}, \t{M}, \t{{\hat{F}}}$ and $\t{K}$ for the
function $\t{f}_t$. At the same time, we also use the notations
in \S\ref{bouquet22} and \S\ref{lemma5}.

Note that, by Lemma \ref{basic lemma},  for $\ve _0>0, \eta_0>0$ and 
$\Delta$ small,
for each $t\in \Delta$, $f_t: {\cB}\cap
f_t^{-1}(\bar{D}_0\setminus\{b_0, \ldots, b_{\sigma}\}) \lar
\bar{D}_0\setminus\{b_0, \ldots, b_{\sigma}\}$ is a
 locally trivial topological fibration.
 Hence
$$\t{f}_t:\t{{\cB}}\cap\t{f}_t^{-1}(\bar{D}_0\setminus\{b_0, \ldots, 
b_{\sigma}\}) \lar \bar{D}_0\setminus\{b_0, \ldots, b_{\sigma}\} $$
is also a locally trivial topological fibration.  
 Denote $\t{E}=\t{{\cB}}\cap \t{f}_t^{-1}(\bar{D}_0)$, $\t{E}^i=\pi 
^{-1}(E^i)$,  $\t{F}^i=\pi ^{-1}(F^i)$ and 
$\t{F}^*=\pi ^{-1}(F^*)$. 
 For $i>0$, $\t{E}^i$ is
a  disjoint union of closed sets, and each of which  is  homeomorphic to
$E^i$. And
$\t{F}^i$ is a  disjoint union of
 closed sets, each of  which is homeomorphic to
$F^i=X\cap \bar{B}^i\cap f_t^{-1}(u_i)$.
These  are possible since the restriction of  $\pi$
to $\t{{\cB}}\setminus * $  is a universal covering.

Then the following proposition is similar to Proposition \ref{additivity}
and can be proved in the same way.

\subsection{Proposition }\label{additivity2}
{\sl
 With the notations and assumptions above,
 we have
$$H_{*}(\t{E}, \tilde{\hat{F}})\cong \bigoplus\limits_{i=0}^{\sigma}
H_{*}(\t{E}^i, \t{F}^i).\eqno{\Box}$$
}

\vs

Note that the ball $B^0$ is a  Milnor ball of $f_t$ at $P_0=0$.
 We have $F^0=F_t$, the Milnor fibre of $f_t$. The lifting
 $\t{F}^0=\t{f}_t^{-1}(u_0)\cap \t{B}^0$ of $F^0$ by $\pi$
 is a universal covering of
 $F^0$ by  the sublemma,
where $\t{B}^0$ is the lifting  of $\bar{B}^0\cap X_t$ by $\pi$.

The fibre $\t{F}^*$ contains
$\t{F}^0$  and  $\t{F}^i \, \,  (i>0)$. For each $i>0$,
$\t{F}^i$ is a disjoint union of pieces which are copies of  $F^i$.
And each  $F^i$ is homotopy equivalent to 
 a bouquet of spheres $F^i\,\, {\buildrel{h}\over{\simeq}} \,\,
S^n\vee \cdots \vee S^n$  ($\beta _i$ copies).
Let $h^i_1\cup \cdots \cup h^i_{\beta _i}$ be the union of the thimbles.
 Denote
$H=\bigcup\limits_{j=1}^{\sigma}(h^j_1\cup \cdots \cup h^j_{\beta _j})$.
Let $x_i$ be the wedge point in $F^i$, and $x_0\in \partial\t{F}^0$.
Let $\tilde{\gamma}_i$ be the union of 
the paths in $\t{F}^*$ connecting $x_0$ and
the liftings of $x_i$ in a usual way. One can take the liftings of $c_i$
 as $\tilde{\gamma}_i$. We have the inclusion:
$$\t{F}':=\t{F}^0\cup(\tilde{\gamma}_1\cup \t{F}^1)\cup\cdots \cup
(\tilde{\gamma}_{\sigma}\cup \t{F}^{\sigma})
\hookrightarrow \t{F}^*.$$
Obviously $\t{F}'$ is homotopy equivalent to $\t{F}^0\vee \t{S}$, where
$\t{S}$ is a wedge of the lifted  bouquets in $F^i$'s ($i>0$).

Denote by  $\t{H}=\bigsqcup H$ the disjoint union of $H$ such that
the attachments of the balls in $\t{H}$ to the spheres in
$\t{S}$ will kill all the the $n$-spheres in $\t{F}^*$ coming from the 
liftings
of $F^i$ ($i>0$). The result of this attachment is denoted by
$\t{F}^*\cup \t{H}$.

Let $B'_j$ be the ball with boundary $S^n_j$, then we have the inclusion
$$F^0\hookrightarrow F^0\vee B':=F^0\vee B_1'\vee \cdots \vee
B'_{\beta},
\quad \text{with } \beta =\sum\limits_{i=1}^{\sigma}\beta _i.$$
Denote by $\t{B}'$  the disjoint union
$\bigsqcup ( B_1'\vee \cdots \vee B'_{\beta})$.
Using the union of  the paths $\tilde{\gamma}_i$  above, we have
$\t{\varphi}: \t{F}^0\vee \t{S} \hookrightarrow \t{F}^*$, the
composition of
$\t{F}^0\vee \t{S} 
\, \,  {\buildrel{h}\over{\simeq}} \,\, 
\t{F}' \hookrightarrow \t{F}^* $,
and 
$\t{F}^0\hookrightarrow \t{F}^0\vee\t{B}'. $ 
We also have the following obvious mappings
$$
\t{\varphi}': \t{F}^0\vee \t{B}' \lar \t{F}^*\cup \t{H}, \quad
\t{\varphi}'':\t{F}^0  \hookrightarrow \t{F}^0\vee \t{B}' \lar \t{F}^*\cup
\t{H}.$$
With the data above, we have the following conclusion  similar to
Theorem \ref{bouquet1} and its proof is  almost word by  word
the same as that of Theorem \ref{bouquet1}.

\subsection{Proposition}\label{covering}
{\sl
The map
$$\t{\varphi}: \t{F}^0\vee \t{S} \lar \t{F}^* $$
induces isomorphisms on all the homology groups:
$$H_*(\t{F}^0\vee \t{S})\cong H_*(\t{F}^*).\eqno{\Box}$$
}

\section{The L\^{e}-Ramanujam problem}\label{section5}

\subsection{}\label{5.1}\hskip -2mm
 Let $\ve _t>0, \eta _t>0$
be admissible for the Milnor fibration of the germ of $f_t$ at 0. Denote
by $B_t$ the open ball in ${\bC}^m$ with center 0 and radius $\ve _t$,
 and by $D_t$  the open disc in ${\bC}$ with center 0 and radius $\eta _t$.
Denote by $\Delta$ a small open disc in ${\bC}$ with center 0.

\subsection{Theorem}\label{le-rama}
{\sl
 Let $(X_t, f_t)$ be a one-parameter deformation of the
isolated singularity $(X_0, f_0)$ with $X_t$ irreducible at $0$ and $\dim X_t\ne 3$.
  Suppose there exists an open  neighborhood U of $0$  such that, for each
$t\in \Delta$, $U\cap X_t\setminus \{0\}$ is non-singular.
  Suppose further that $(X_t ,0)$ has link stability and the 
 isolated singularity $(X_0,f_0)$  does not split. Then
\begin{itemize}
\item[1)]{the homotopy type of the Milnor fibre of $f_t$ is constant;
  i.e.  for any $u\in \bar{D}^*_t :=\bar{D}_t\setminus \{0\}$
$$F_0=
f_0^{-1}(u)\cap X_0\cap \bar{B}_0
{\buildrel{h}\over{\simeq}}
F_t=f_t^{-1}(u)\cap X_t\cap \bar{B}_t; $$}
\item[2)]{The   monodromy fibrations  of $f_0$ and $f_t$
  are fibre  homotopy equivalent (as fibrations over
$\partial\bar{D}_0$ and $\partial\bar{D}_t$ respectively); i.e.
$$E_0:=f_0^{-1}(\partial\bar{D}_0)\cap X_0\cap \bar{B}_0
\quad {\buildrel{h}\over{\simeq}}\quad
E_t:=f_t^{-1}(\partial\bar{D}_t)\cap X_t\cap \bar{B}_t.$$
   Here and in the following, we
 denote the fibration $(E_t,   f_t|{E_t},\partial\bar{D}_t)$ by $E_t$;}
\end{itemize}

If, moreover,  the Milnor fibre $F_t$ of $f_t$ is simply connected or if $\dim X_t = 2$, 
then we have:

\begin{itemize}
\item[3)]{ the  diffeomorphism type of
 the Milnor fibration of $f_t$ is constant (as fibrations over
$\partial\bar{D}_0$ and $\partial\bar{D}_t$ respectively); i.e.
$$E_0:=f_0^{-1}(\partial\bar{D}_0)\cap X_0\cap \bar{B}_0
\quad {\buildrel{\scriptstyle\text{diffeo}}\over{\simeq}}\quad
E_t:=f_t^{-1}(\partial\bar{D}_t)\cap X_t\cap \bar{B}_t;$$
 }
\item[4)]{the local topological type of $f_t$ is constant; i.e.,
$$\left(\bar{B}_0\cap X_0, \bar{B}_0\cap X_0\cap f_0^{-1}(0)\right)
{\buildrel{\scriptstyle\text{homeo}}\over{\simeq}}
\left(\bar{B}_t\cap X_t, \bar{B}_t\cap X_t\cap f_t^{-1}(0)\right).$$
}
\end{itemize}

}

\prf
 We follow the pattern of  L\^{e}-Ramanujam's proof \cite{LR}.
By the non-splitting condition,
$f_t$ has no critical point on $U\cap X_t\setminus \{0\}$
 for any $t$.

 Note that the diffeomorphism type of
the Milnor fibration does not depend on
the choice of  $\eta>0$.
 Hence for any $t\ne 0$ fixed, as fibrations over $\partial\bar{D}_0$
and $\partial\bar{D}_t$ respectively, we have a diffeomorphism
$$E_0=f_0^{-1}(\partial\bar{D}_0)\cap X_0\cap\bar{B}_0
\quad {\buildrel{\scriptstyle\text{diffeo}}\over{\simeq}}\quad
E'_0:=f_0^{-1}(\partial\bar{D}_t)\cap X_0\cap\bar{B}_0.$$
Then we prove that (as fibrations)
$$E_0'=
f_0^{-1}(\partial\bar{D}_t)\cap X_0\cap \bar{B}_0
\quad {\buildrel{\scriptstyle\text{diffeo}}\over{\simeq}}\quad
E'_t:=f_t^{-1}(\partial\bar{D}_t)\cap X_t\cap \bar{B}_0$$
So there is an inclusion $E_t \hookrightarrow E'_t$, and also
inclusions for their fibres.
We prove these induce the desired results.

Consider  the map $G$ defined in the proof of Lemma \ref{basic lemma}
$$G(z, t)=(f(z, t), t):X\cap(\bar{B}_0\times\Delta)\longrightarrow {\bC}
\times \Delta.$$
This map induces the following two
differentiable fibrations by Ehresmann's theorem (see \cite{La}):
$$G_1:\left(\bar{B}_0\times \Delta\right)\cap X\cap
G^{-1}\left(\partial\bar{D}_t\times \Delta \right)\longrightarrow
\partial\bar{D}_t\times \Delta,$$
and
$$G_2:\left(\partial\bar{B}_0\times \Delta\right)\cap X\cap
G^{-1}\left(\bar{D}_t\times \Delta \right)\longrightarrow
\bar{D}_t\times \Delta.$$
Moreover, $G_2$ is a trivial fibration since
 $\bar{D}_t\times \Delta$ is contractible.

Hence, as fibrations over $ \partial\bar{D}_t\times 0$ and
$\partial\bar{D}_t\times t$ respectively,
$$G^{-1}(\partial\bar{D}_t\times 0)
\quad {\buildrel{\scriptstyle\text{diffeo}}\over{\simeq}}\quad
G^{-1}(\partial\bar{D}_t\times t),$$
and this is compatible with the trivialization $G_2$. This proves that
as fibrations over $\partial\bar{D}_t$
$$E_0'=
f_0^{-1}(\partial\bar{D}_t)\cap X_0\cap \bar{B}_0
\quad {\buildrel{\scriptstyle\text{diffeo}}\over{\simeq}}\quad
E'_t:=f_t^{-1}(\partial\bar{D}_t)\cap X_t\cap \bar{B}_0.$$

Next, we prove that
$E'_t
\quad {\buildrel{\scriptstyle\text{diffeo}}\over{\simeq}}\quad E_t$
  as fibrations over $ \partial\bar{D}_t$.
Recall that  a fibre of $E_t$ is denoted by $F_t$, which is the
 Milnor fibre of $f_t$.
Since the fibre of $E'_t$ is diffeomorphic  to the
 Milnor fibre $F_0$ of $f_0$, in the following we use this notation.


Obviously
$E_t\hookrightarrow E'_t$.
By Ehresmann's theorem (loc. cit.)
$$f_t: \triangledown E_t:=\overline{E'_t\setminus E_t}
\longrightarrow\bar{D}_t $$
is a differentiable fibration and trivial:
$\triangledown E _t\quad 
{\buildrel{\scriptstyle\text{diffeo}}\over{\simeq}}\quad
\text{fibre}\times \bar{D}_t,$ as fibrations over $ \bar{D}_t.$
 We may take the typical  fibre to be the one  over  $u \in \partial\bar{D}_t 
$:
$$\triangledown F_t  :=\triangledown E_t\cap f^{-1}_t(u)=F_0\setminus
 \text{int}F_t. $$

We use the following lemma  which will be proved in
\S\ref{ProofofLemma5.1}.

\subsection{Lemma}\label{lemma5.1}
{\sl Under the assumptions of Theorem \ref{le-rama} and with the notations above, we have:
\begin{itemize}
\item[a)] { If $\dim X_t >3$, then the inclusion $F_t\hookrightarrow F_0$
 is a homotopy equivalence  (cf. \cite{S2});}
\item[b)] {If $\dim X_t >3$, then $\pi _1(\partial F_0)\cong\pi _1(\triangledown F_t)\cong
 \pi _1(\partial F_t)$   (loc. cit.);}
\item[c)]{$
H_{*}(\partial F_t, {\bZ})\cong H_{*}(\triangledown F_t, {\bZ}) $;}
\item[d)]{$
H_{*}(\partial F_0, {\bZ})
\cong H_{*}(\triangledown F_t, {\bZ})$.}
\end{itemize}
All the isomorphisms are induced by appropriate inclusion maps.
}

\vs

In case $\dim X_t=2$, the statements 1) and 2) follow from 
Theorem~\ref{bouquet22}
and the fact that the Milnor
fibres are bouquets of 1-spheres.  

We now consider the case  $\dim X_t>3$.
The statement 1) follows from a) in   Lemma \ref{lemma5.1},  and 2) follows
 from  1) and a theorem of Dold  \cite[(6.3)]{dold}.

3) By Morse theory,   $\pi _1(\partial F_0)=\pi _1(F_0) $ and 
$\pi _1(\partial F _t)=\pi _1(F_t) $. Hence $\partial F_0$ and
$\partial F_t$ are simply connected
 since  $F_0$ and $F_t$ are simply connected by the assumption. By Whitehead's theorem
$\partial F_0 \hookrightarrow \triangledown F_t$ and
 $\partial F _t\hookrightarrow \triangledown F_t$ are homotopy equivalences.
Since $\dim _{\bR} \triangledown F_t\ge 6$, by $h$-cobordism theorem (cf.
\cite{Milnor2, smale}),
$\triangledown F_t
\quad {\buildrel{\scriptstyle\text{diffeo}}\over{\simeq}}\quad
[0,1]\times \partial F_t.$
Since
$$f_t:\partial\bar{B}_0\cap X_t\cap f^{-1}_t (\bar{D}_t)
\longrightarrow \bar{D}_t$$
is trivial,
$$f_t:\partial\bar{B}_0\cap X_t\cap f^{-1}_t (\partial\bar{D}_t)
\longrightarrow \partial\bar{D}_t$$
is also trivial.
Hence $E'_t$ can be  obtained from $E_t$ by attaching
 $\partial\bar{D}_t\times \text{collar}$. This proves that
 $E'_t
\quad {\buildrel{\scriptstyle\text{diffeo}}\over{\simeq}}\quad E_t$
as fibrations over $\partial\bar{D}_t$.

4) Let $\Phi_t :  E'_t \longrightarrow E_t $
and  $\lambda_t:[0,1]\times \partial F_t\times \bar{D}_t
\longrightarrow \triangledown E_t=\triangledown F_t\times\bar{D}_t $
be the diffeomorphisms of fibrations obtained
 in 3). Assume $\Phi_t (\lambda_t (0,z,u))=\lambda_t (1,z,u)$.
We also have a diffeomorphism of fibrations
$$\Psi_t :\partial \bar{B}_0\cap X_t\cap f^{-1}_t(\bar{D}_t)
\longrightarrow
\partial \bar{B}_t\cap X_t\cap f^{-1}_t(\bar{D}_t)$$
and  $\Phi_t $ and $\Psi_t$ are
equal at the points where both of them are defined. Furthermore
$$\lambda_t (0\times \partial F_t\times 0)
= \partial \bar{B}_0\cap X_t\cap f^{-1}_t(0)
\quad {\buildrel{\scriptstyle\text{diffeo}}\over{\simeq}}\quad
\lambda_t (1\times \partial F_t\times 0)
= \partial \bar{B}_t\cap X_t\cap f^{-1}_t(0)$$
under $\Psi_t$.
By \cite[Proposition 5.4]{lo} and its proof, there is a homeomorphism
$$
\left[\bar{B}_t\cap X_t\cap f^{-1}_t(\partial\bar{D}_t)\right]
\cup \left[ \partial \bar{B}_t\cap X_t\cap f^{-1}_t(\bar{D}_t)\right]
\longrightarrow
\partial\bar{B}_t\cap X_t$$
preserving  $ \partial\bar{B}_t\cap f^{-1}_t(0) $.
We have 
$$\partial\bar{B}_t\cap X_t\cap f^{-1}_t(0)
{\buildrel{\scriptstyle\text{homeo}}\over{\simeq}}
\partial\bar{B}_0\cap X_t\cap f^{-1}_t(0)
{\buildrel{\scriptstyle\text{homeo}}\over{\simeq}}
\partial\bar{B}_0\cap X_0\cap f_0^{-1}(0),$$
where  the second homeomorphism comes from $G_2$.
Hence
$$\begin{array}{ll}
& \left(\bar{B}_t\cap X_t, \bar{B}_t\cap X_t\cap f^{-1}_t(0)\right)\cr
{\buildrel{\scriptstyle\text{homeo}}\over{\simeq}}&
\left(\bar{B}_t\cap X_t,
\text{Cone}(\partial\bar{B}_t\cap X_t\cap f^{-1}_t(0))\right)\cr
{\buildrel{\scriptstyle\text{homeo}}\over{\simeq}}&
\left(\bar{B}_0\cap X_t,
\text{Cone}(\partial\bar{B}_0\cap X_t\cap f^{-1}_t(0))\right)\cr
{\buildrel{\scriptstyle\text{homeo}}\over{\simeq}}&
\left(\bar{B}_0\cap X_0,
 \text{Cone}(\partial\bar{B}_0\cap X_0\cap f_0^{-1}(0))\right)\cr
{\buildrel{\scriptstyle\text{homeo}}\over{\simeq}}&
\left(\bar{B}_0\cap X_0,
 \bar{B}_0\cap X_0\cap f_0^{-1}(0)\right),\cr
\end{array}$$
where the first and the  last  homeomorphisms were proved by Iomdin
 \cite{Iomdin}, the second and the third homeomorphisms follow from the
 above discussions.
 
 Finally, returning to the case $\dim X_t = 2$, the proof of 3) and 4) follows from c) and d) of Lemma \ref{lemma5.1} together with the fact that a real two-dimensional homology cobordism is a product.
\eprf

\subsection{Proof of Lemma \ref{lemma5.1}}\label{ProofofLemma5.1}

The proof of  b) is  essentially contained
 in the proof of Lemma \ref{lemma5}.

a) In this special case, the map $\varphi $ defined in \S\ref{2.3}
is in fact the inclusion map  $F_t\hookrightarrow F_0$.
It follows from Lemma \ref{lemma5}  that
the inclusion $F_t\hookrightarrow F_0$ induces an isomorphism
of their fundamental groups.  By  Proposition \ref{covering},
it also induces an isomorphism on the homology of the universal coverings
of the spaces.
Then we use  Whitehead's theorem \cite{Wh}.

c) Since  $H_{*}(F_0, F_t) $ is trivial, by using
 excision theorem we have   $H_{*}(\triangledown F_t, \partial F_t)=0$. The
proof of c)
 is finished.

 d) The isomorphism in d) mainly comes from the Poincar\'e-Lefschetz duality
theorem (see e.g. \cite{Span}).
\hspace*{\fill}$\Box$
 \subsection{Example}\label{5.3}  Let $X_t$ be the  Brian\c{c}on-Speder
surfaces
(\cite{BS}, see also Example \ref{example4-3}) defined by
$h_t=z^3+ty^{2\alpha+1}z+xy^{3\alpha +1}+x^{6\alpha +3}=0$ $(\alpha \ge 1)$.
  These  surfaces are quasi-homogeneous.
Consider the function
$f_t=xy^{\alpha}+z+tz^2$
on $X_t$. The  critical locus of $f_t$ is the solutions of the following
system of equations:
\begin{align}
\alpha (3z^2+ty^{2\alpha+1})xy^{\alpha -1}-\left((2\alpha+1)ty^{2\alpha}z+
(3\alpha+1)xy^{3\alpha}\right)(1+2tz)=0\tag{1}\\
(3z^2+ty^{2\alpha +1})y^{\alpha}-\left((6\alpha +3)x^{6\alpha+2}+y^{3\alpha+1}
\right)(1+2tz)=0\tag{2}\\
ty^{3\alpha}z+xy^{4\alpha}-3\alpha x^{6\alpha +3}y^{\alpha -1}=0\tag{3}\\
xy^{\alpha}+z+2tz^2=0\tag{4}\\
h_t=0\tag{5}
 \end{align}
The equations  (1)--(3) come from the minors of the Jacobian of $(f_t,h_t)$
and  (4) comes from the differential of $f_t$ by the Euler
derivation.

 Note that if any one of the $x, y $, or $z$ is zero, then the other two are
also zero. So we may assume $xyz\ne 0$.
From (3), (4) and (5), one has
$$z=\omega x^{2\alpha +1},\quad y^{\alpha}=-\omega ux^{2\alpha},\quad 
3\alpha x^2=\omega ^3u^2(t-u)y,$$
where $u:=1+2tz $ and $\omega ^3:=-(3\alpha+1).$
From this we obtain that
$tz=c(t)$ with  $c(0)\ne 0$. This means that the non-zero solutions of the 
 above equations tend
to infinity as $t$ tends to 0. We conclude that, for small $t$,
the singularity of $(X_0,f_0)$ does not split in a neighborhood of the origin.
By Theorem \ref{le-rama},
 the local topological type of $f_t$ and 
the homotopy type of
 the Milnor fibration of $f_t$ are constant.

\end{document}